\documentclass[conference]{IEEEtran}

\ifCLASSINFOpdf

\else

\fi

\usepackage{graphicx}
\usepackage{amsmath,amssymb,epsfig,epstopdf}

\begin{document}
\begin{titlepage}\centering
		\vspace*{\fill}
		\centering
		\LARGE \begingroup
		Pre-print version \\ of \\ 
		"RVP-FLMS : A Robust Variable Power Fractional LMS Algorithm" \\by \\
		Jawwad Ahmad, Muhammad Usman, Shujaat Khan, Imran Naseem and Hassan Jamil Syed \\
		submitted to \\
		IEEE International Conference on Control System, Computing and
		Engineering (ICCSCE). IEEE, 2016.
		\endgroup
		\vspace*{\fill}
\end{titlepage}

\title{RVP-FLMS : A Robust Variable Power Fractional LMS Algorithm}


\author{\IEEEauthorblockN{Jawwad Ahmad\IEEEauthorrefmark{1},
		Muhammad Usman\IEEEauthorrefmark{1}, Shujaat Khan\IEEEauthorrefmark{1}\IEEEauthorrefmark{2}, Imran Naseem\IEEEauthorrefmark{3} and
		Hassan Jamil Syed\IEEEauthorrefmark{4}}
	\IEEEauthorblockA{\\
		\IEEEauthorrefmark{1}Faculty of Engineering Science and Technology (FEST),\\ Iqra University, Defence View,\\ Karachi-75500, Pakistan.\\
		Email: \{jawwad, musman, shujaat\}@iqra.edu.pk\\
		\\
		\IEEEauthorrefmark{2}Department of Bio and Brain Engineering,\\ Korea Advanced Institute of Science and Technology (KAIST),\\ 335 Gwahangno, Yuseong-gu, Daejeon 305-701, Korea.\\
		Email: shujaat@kaist.ac.kr\\
		\\
		\IEEEauthorrefmark{3}College of Engineering,\\ Karachi Institute of Economics and Technology,\\ Korangi Creek, Karachi 75190, Pakistan.\\
		Email: imrannaseem@pafkiet.edu.pk \\ 
		\\
		\IEEEauthorrefmark{4}Faculty of Computer Science and Information Technology,\\ University of Malaya, Malaysia.\\
		Email: shjamil@siswa.um.edu.my}}
\maketitle
%\author{
%\IEEEauthorblockN{Jawwad Ahmad, Muhammad Usman} % and Shujaat Khan}
%\IEEEauthorblockA{Faculty of Engineering Science and Technology,\\
%Iqra University, Defence View, Karachi-75500, Pakistan.\\
%Email: jawwad@iqra.edu.pk, musman@iqra.edu.pk}%, 
%%\\shujaat@iqra.edu.pk}
%\and
%\IEEEauthorblockN{Shujaat Khan}
%\IEEEauthorblockA{Department of Bio and Brain Engineering, \\
%Korea Advanced Institute of Science \& Technology (KAIST),\\
%335 Gwahangno, Yuseong-gu, Daejeon 305-701, Korea. \\
%Email: shujaat@kaist.ac.kr}
%\and
%\IEEEauthorblockN{Imran Naseem}
%\IEEEauthorblockA{College of Engineering (CoE), \\
%Karachi Institute of Economics and Technology, \\
%Korangi Creek, Karachi 75190, Pakistan.\\
%Email: imrannaseem@pafkiet.edu.pk}
%\and
%\IEEEauthorblockN{Hassan Jamil Syed \IEEEmembership{Senior~Member,~IEEE}}
%\IEEEauthorblockA{Faculty of Computer Science and Information Technology,\\University of Malaya, Malaysia. \\
%Email: shjamil@siswa.um.edu.my}
%}
%% make the title area
%\maketitle

\begin{abstract}
In this paper, we propose an adaptive framework for the variable power of the fractional least mean square (FLMS) algorithm.  The proposed algorithm  named as robust variable power FLMS (RVP-FLMS) dynamically adapts the fractional power of the FLMS to achieve high convergence rate with low steady state error.  For the evaluation purpose, the problems of system identification and channel equalization are considered.  The experiments clearly show that the proposed approach achieves better convergence rate and lower steady-state error compared to the FLMS.  The MATLAB code for the related simulation is available online at https://goo.gl/dGTGmP.\\
\end{abstract}
\begin{IEEEkeywords}
	\normalfont{
Least mean square (LMS), fractional calculus, plant identification, channel equalization, fractional LMS (FLMS), robust variable power FLMS (RVP-FLMS), robust variable step size (RVSS), high convergence, low steady state error, adaptive filter.}
\end{IEEEkeywords}

\IEEEpeerreviewmaketitle

\section{Introduction}
The conventional integer value integration and differentiation is undoubtedly a fundamental tool of calculus being utilized by professionals to determine the solution of numerous problems in various fields.  On the contrary fractional order calculus (FOC) deals with the derivative and integral calculated with the fractional power.
%  The concept of FOC is first defined by Leibniz, the notation of FOC that he used in his paper is $\frac{D^{n}x}{Dx^{n}}$ where $n=1/2$ \cite{FR1}.  
Derivative with the fractional power is certainly a paradox which introduced a new aspect in mathematics.  Though the concepts of FOC were developed earlier in the $16^{th}$ century by Liouville, Reimann and Leibniz \cite{liou}, tendency of its utilization in generic areas became notable in the past decade.  With the maturation of theoretical and practical operators of FOC,  researchers have shown an inclination towards the implementation of FOC to observe the solution of the problems that were solved using the traditional calculus \cite{FCA6}. 

For the problem of path tracking in autonomous vehicles FOC was used in \cite{FCA1}. FOC was used  in \cite{FCA2} to give an experimental proof of a theoretical model proposed for the pulse propagation through porous medium.  FOC has turned out to be a convenient tool in the applications of viscoelasticity \cite{FCA3}, image processing \cite{FCA4}, edge detection \cite{FCAap}, control systems and many others \cite{FCAap2}.  Implementation of FOC in the field of signal processing has been suggested and opted by a number of researchers.  In this context, the utilization of FOC in least mean square (LMS) adaptive filtering is being particularly focused.

FOC based least mean square algorithm was proposed in \cite{LMS12} as FLMS.  The algorithm was implemented on system identification problem and the results were compared with LMS.  Several variants of FLMS algorithm have been proposed showing better performance \cite{LMS14, LMS16, LMS17}.  In this paper a contemporary modification of the FLMS algorithm is proposed.  The proposal is to use a variable fractional power rather than a fixed one, the concept of which is inspired by variable learning rate least mean square algorithm (RVSS-LMS) \cite{RVSS1}.  The objective is to use variable power of the fractional derivative.  The proposed algorithm is evaluated on the problem of equalization and plant identification, the results of which are compared with the conventional LMS.  The rest of the paper is structured as follows: A brief literature review of modified least mean algorithms is discussed in section \ref{FLMS}. The details of proposed algorithm is discussed in section \ref{RVFPflms} followed by the experimental setup in section \ref{Experimental}. The conclusion of the paper is drawn in section \ref{Conclusion}. 

\section{Least Mean Square}\label{FLMS}
Least mean square algorithm is an adaptive filter that belongs to the group of stochastic gradient methods.  LMS has been extensively adopted to counter a number of problems in the area of signal processing.  A considerable amount of work is done towards the mitigation of the drawbacks in LMS.  

Researchers have used several techniques to optimize the performance of LMS \cite{LMS14, LMS13, LMS11, LMS9, LMS15}.  The theory of FOC was applied to the LMS algorithm and a fractional least mean square (FLMS) was proposed in \cite{LMS15}.  The algorithm was tested for the estimation of input nonlinear control autoregressive (INCAR) models and the output resulted in better convergence rates when compared to volterra least mean square (VLMS) and kernel least mean square (KLMS).  The modified fractional least mean square (MFLMS) was proposed in \cite{LMS16}.  The algorithm was designed to improve the computational complexity involved in FLMS due to the complex gamma function.  MFLMS has therefore shown improved results compared to the LMS and FLMS. In another work \cite{LMS17} the adaptive weight gain parameters are incorporated.  Gradient-based approach is implemented on variable learning scheme that changed the nature of the order of fractional derivative in MFLMS from fixed to adaptive.  They proposed adaptive step-size modified fractional least mean square (AMFLMS) and presented improvement by comparing their algorithm with the LMS, FLMS and MFLMS.  

In this research we propose to utilize the concept of robust variable step size (RVSS) \cite{RVSS1} for the variable fractional power of FLMS algorithm.  The proposed scheme is robust and computationally less expensive.  For the performance evaluation we consider the problem of system identification and equalization. The results are compared with the FLMS algorithm.

\section{Proposed RVP FLMS}\label{RVFPflms}
In \cite{LMS12} the weight update equation for FLMS is given as :
	
	\begin{eqnarray}\label{weightupdate1}
	w_{k}(n+1)= w_{k}(n)- \mu \dfrac{\partial J(n)}{\partial w_{k}} - \mu_{f}\left(\dfrac{\partial}{\partial w_{k}}\right)^v J(n)
	\end{eqnarray}
	
	where $w_{k}(n)$ is the weight of the $k^{th}$ tap at $n^{th}$ iteration, $J(n)$ is the cost function defined in equation (\ref{cost}), $v$ is the fractional power of derivative, $\mu$ and $\mu_{f}$ are the step sizes.
	
	\begin{eqnarray}\label{cost}
	J(n)= \dfrac{1}{2}e(n)^2 = \dfrac{1}{2}(d(n)-y(n))^{2}
	\end{eqnarray}

where $e(n)$ is the instantaneous error between the desired output $d(n)$ and the estimated output $y(n)$ at $n^{th}$ iteration.
	
According to chain rule $\dfrac{\partial J(n)}{\partial w_{k}}$ and $\left(\dfrac{\partial}{\partial w_{k}}\right)^v J(n)$ are defined as :

	\begin{eqnarray}\label{costd}
	\dfrac{\partial J(n)}{\partial w_{k}} = \dfrac{\partial J(n)}{\partial e} \dfrac{\partial e(n)}{\partial y} \dfrac{\partial y(n)}{\partial w_{k}}
	\end{eqnarray}
	
	\begin{eqnarray}\label{costdf}
	\left(\dfrac{\partial}{\partial w_{k}}\right)^v J(n) = \dfrac{\partial J(n)}{\partial e} \dfrac{\partial e(n)}{\partial y} \left(\dfrac{\partial}{\partial w_{k}}\right)^v y(n)
	\end{eqnarray}

solving equation \eqref{costd} result in equation (\ref{costd1})
	\begin{eqnarray}\label{costd1}
	\dfrac{\partial J(n)}{\partial w_{k}} = -e(n)x(n)
	\end{eqnarray}

Using the Rieman-Lioville fractional derivative method
\begin{eqnarray}
(D^v f)(t)= \dfrac{1}{\Gamma}(\dfrac{d}{dt})^n \int\limits_{0}^{t}(t-\tau)^n-v-1 f(\tau)d\tau
\end{eqnarray}

\begin{eqnarray}
D^v(t-a)^\alpha = \dfrac{\Gamma(1+\alpha)}{\Gamma(1+\alpha+v)}(t-a)^\alpha-v
\end{eqnarray}

equation \eqref{costdf} result in equation (\ref{costdf1})

	\begin{eqnarray}\label{costdf1}
	\left(\dfrac{\partial}{\partial w_{k}}\right)^v J(n) = -e(n)x(n)\dfrac{w_{k}^{(1-v)}(n)}{\Gamma(2-v)}
	\end{eqnarray}	 
	
Using equation \eqref{costd1} and  \eqref{costdf1} the weight update equation \eqref{weightupdate1} becomes : 
	
	\begin{multline}\label{weightupdate2}
	w_{k}(n+1)=w_{k}(n) + \mu e(n)x(n) \\ + \mu_{f}e(n)x(n)\dfrac{w_{k}^{1-v}(n)}{\Gamma(2-v)}
	\end{multline}
%	we assume $w_{k}^{1-f}(n)\cong w_{k}^{1-f}(n-1)$
%
%	
%	\begin{multline}
%		w_{k}(n)=w_{k}(n-1) + \mu e(n)x(n-k) \\ + \mu_{f}e(n)x(n-k)\dfrac{w_{k}^{1-f}(n-1)}{\Gamma(2-f)}
%	\end{multline}

%%%	For simplicity, in equation \eqref{weightupdate2}, we consider $\mu_{f} = \mu \Gamma(2-v)$ and this is result in :
%%	
%%		\begin{multline}
%%			w_{k}(n+1)=w_{k}(n) + \mu e(n)x(n) \\ + \mu e(n)x(n) w_{k}^{1-v}(n)
%%		\end{multline}
%%	
%		\begin{equation}
%			w_{k}(n+1)=w_{k}(n) + \mu e(n)x(n) \left(1 + w_{k}^{1-v}(n)\right)
%		\end{equation}
%	
	For time varying fractional power the $v$ can be replaced by $v(n)$ :
	\begin{multline}\label{weightupdate2}
	w_{k}(n+1)=w_{k}(n) + \mu e(n)x(n) \\ + \mu_{f}e(n)x(n)\dfrac{w_{k}^{1-v(n)}(n)}{\Gamma(2-v(n))}
	\end{multline}

The update rule for the time varying fractional power $v$ using the RVSS-LMS based method is defined as \cite{RVSS1} :	
		\begin{eqnarray}
		v (n+1)= \beta v(n)+ \gamma p^2(n)
		\end{eqnarray}	

	where $(0 < \beta < 1)$, $(\gamma > 0)$, $p(n)$ is the average error energy correlation and $v(n+1)$ is set to $v_{min}$ or $v_{max}$
	
	\begin{equation}
	p(n) = \alpha p(n - 1) + (1 - \alpha) e(n) e(n - 1)
	\end{equation}
	
	The positive constant $ \alpha \; (0 < \alpha < 1) $ is a weighting parameter that governs the averaging time constant. It is in actual the forgetting factor.  The limits on $ v (n + 1) $ is given by:
	\begin{equation}
	v (n + 1) = \left\{ \begin{array}{rcl}
	v_{\max} & \mbox{if} & v (n + 1) > v_{\max} \\ 
	v_{\min} & \mbox{if} & v (n + 1) < v_{\min} \\ 
	v(n + 1) & & otherwise
	\end{array}\right.
	\end{equation}
	
	where $ v_{\max} > v_{\min} > 0 $. 
	
	\section{Simulation Setup and Results}\label{Experimental}
	
	 For the evaluation of the proposed method the problem of system identification and equalization are considered. 
\subsection{System Identification}\label{plantidentification}

	 To estimate the mathematical model of a system, the system identification method is used. Literature including \cite{LMS12,LMS13,LMS15,LMS1, khan2016novel} show that the adaptive techniques produce good performance in the application of system identification.  To assess the performance of the proposed RVP-FLMS, we consider a linear system, shown in Figure \ref{plant1}:	
	\begin{figure}[b!]
			\begin{center}
				\centering %\fbox{
				\includegraphics*[scale=0.32,bb=0 0 720 330]{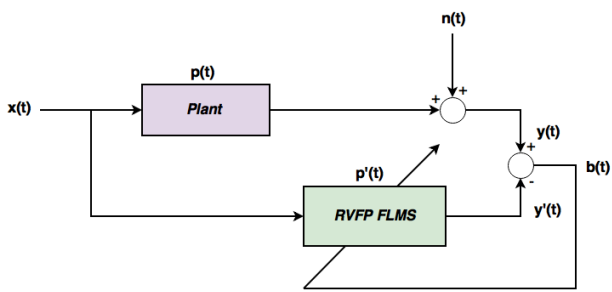}% }
			\end{center}
			\caption{System identification using adaptive learning algorithm.}
			\label{plant1}
		\end{figure}

\begin{multline}\label{planteq}
		y(t)=a_1x(t)+a_2x(t-1)+a_3x(t-2) + n(t)
		\end{multline}
	
Equation (\ref{planteq}) shows the mathematical model of the system, where $x(t)$ is the input and $y(t)$ is the output of the system, the disturbance model $n(t)$ is assumed to be $\mathcal{N}(0,\sigma_d^2)$, $a_i$s  depicts the  polynomial coefficients representing number of zeros in the system.  For the experiment, $x(t)$ is taken to be a binary phase shift keying (BPSK) modulated 500 randomly generated samples.  In Figure \ref{plant1}, the impulse response of the system is $p(t)$ while $\hat{y}(t)$ is approximate output, $\hat{p}(t)$ is the approximate impulse response and $b(t)$ is the error of estimation.  The simulation parameters selected for the experiments are:  $a_1=0.9$, $a_2=0.3$, $a_3=-0.1$.  The experiments are performed on two noise levels with the SNR values of $10$ dB and $20$ dB. 

			\begin{figure}[h!]
				\begin{center}
					\centering%\fbox{
					\includegraphics*[scale=0.4,bb=50 185 560 590]{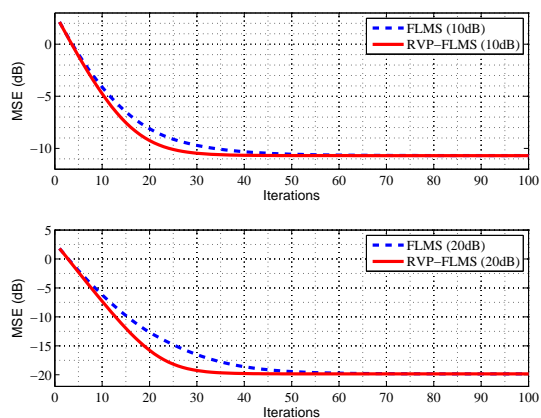}% }
				\end{center}
				\caption{Mean square error (MSE) curve for system identification problem}
				\label{planterror}
			\end{figure}
		
	There are $3$ number of taps for the model structure.  The weights values were initialised to zero. For the FLMS and RVP-FLMS value of the $\mu$ and $\mu_{f}$ was set to $1\times10^{-4}$ and the initial value of the fractional power $v$ is taken to be $0.5$.  For the proposed RVP-FLMS the values of $\alpha$, $\beta$ and $\gamma$ are all set to $0.9$, $0.99$, and $0.9$ respectively. The value of $v_{max}$ and $v_{min}$ are chosen to be $0.5$ and $1$ respectively.
	
	The experiment are performed on two noise levels the mean square error (MSE) curves are shown in Figure \ref{planterror}.  For the SNR values of $10$ dB and $20$ dB the RVP-FLMS shows the best performance. The MSE of $-10.44$ dB and $-20.20$ dB is achieved in around 45 and 50 iterations respectively.  Whereas the FLMS converges to the MSE value of $-10.43$ dB and $-20.17$ dB in around 80 and 90 iterations respectively.
	
	To compare the time complexity of the suggested method with FLMS, we investigated the training time for 100 iterations.  The proposed method utilizes 1.85 seconds whereas the FLMS takes 1.57 seconds.  The experiment clearly shows that the proposed robust approach dynamically adapts the learning rate to achieve the minimum steady state error in lesser number of iteration.

\subsection{Equalization}
	 The equalizer is used to nullify the effect of noise and distortion caused by the channel.  Adaptive learning techniques show good results in this context \cite{EQ1,EQ2,EQ3,EQ4}.  To evaluate the effectiveness of the proposed RVP-FLMS, we consider a linear filter, shown in Figure \ref{EQL}:
	
	\begin{figure}[h!]
			\begin{center}
				\centering %\fbox{
				\includegraphics*[scale=0.45,bb=0 0 460 185]{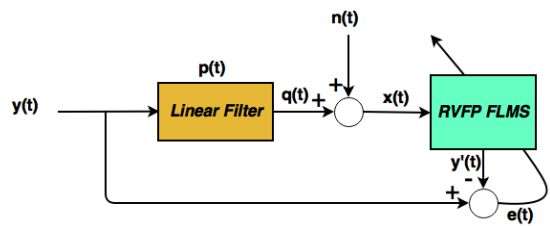} %}
			\end{center}
			\caption{Equalization using adaptive learning algorithm.}
			\label{EQL}
		\end{figure}

\begin{multline}\label{EQLE}
		q(t)=a_1y(t)+a_2y(t-1)+a_3y(t-2)
		\end{multline}
	
Equation (\ref{EQLE}) shows the mathematical model of the linear filter, where $y(t)$ is the input and $q(t)$ is the  output of the filter. $a_i$s depicts the polynomial coefficients which are actually the number of zeros in the system.  For the purpose of this experiment the same plant defined in section (\ref{plantidentification}) is used. 

 In Figure \ref{EQL}, the impulse response of the filter is $p(t)$ while the disturbance model $n(t)$ is supposed to be $\mathcal{N}(0,\sigma_d^2)$.  $\hat{y}(t)$ is the approximate input, $x(t)$ is the approximate noisy output and $e(t)$ is the  error of estimation.  The experiments are performed on two noise levels with the SNR values of $10$ dB and $20$ dB. 

					\begin{figure}[h!]
						\begin{center}
							\centering%\fbox{
							\includegraphics*[scale=0.4,bb=50 185 560 590]{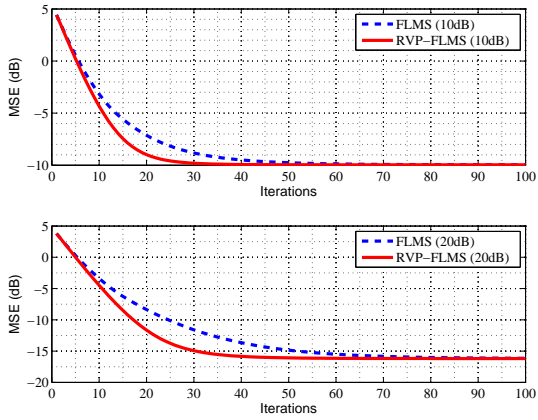} %}
						\end{center}
						\caption{Mean square error (MSE) curve for equalization problem}
						\label{EQLError}
					\end{figure}
		
	In the equalizer structure, $3$ number of taps were chosen.  The weights values were initialised to zero. For the FLMS and RVP-FLMS value of the $\mu$ and $\mu_{f}$ was set to $1\times10^{-4}$ and the initial value of the fractional power $v$ is taken to be $0.5$.  For the proposed RVP-FLMS the values of $\alpha$, $\beta$ and $\gamma$ are all set to $0.9$, $0.99$, and $0.9$ respectively. The value of $v_{max}$ and $v_{min}$ are chosen to be $0.5$ and $1$ respectively.
	
	The experiment are performed on two noise levels the mean square error (MSE) curves are shown in Figure \ref{EQLError}.  For the SNR values of $10$ dB and $20$ dB the RVP-FLMS shows the best performance. The MSE of $-9.96$ dB and $-16.33$ dB is achieved in around 40 and 50 iterations respectively.  Whereas the FLMS converges to the MSE value of $-9.58$ dB and $-16.28$ dB in around 70 and 90 iterations respectively.
	
	To compare the time complexity of the suggested method with FLMS, we investigated the training time for 100 iterations.  The suggested method requires 1.92 seconds whereas the FLMS takes 1.63 seconds.  The experiment clearly shows that the proposed robust approach dynamically adapts the learning rate to achieve the minimum steady state error in lesser number of iteration.
			
	\section{Conclusion}\label{Conclusion}
	In this research an adaptive least mean square algorithm based on fractional derivative is proposed. In particular, the fractional power of the FLMS is made variable using the concept of robust variable step size.   The performance of the proposed approach is evaluated on system identification and equalization problems.  The results of the proposed algorithm are compared to the FLMS.  The proposed algorithm attains better convergence rate and steady-state error and is therefore found to be superior to the FLMS algorithm.

%	% use section* for acknowledgment
	\section*{Acknowledgment}
	The research is partially supported by the National ICT R\&D Fund (http://www.ictrdf.org.pk)
\bibliographystyle{IEEEtran}
\bibliography{bare_conf}

\end{document}